\documentclass[11pt]{article}

\usepackage{latexsym}
\usepackage{amsmath}
\usepackage{amsthm}
\usepackage{amssymb}
\usepackage{vmargin}
\usepackage{amscd}
\usepackage{stmaryrd}
\usepackage{euscript}
\usepackage{mathrsfs}
\usepackage{amscd}
\usepackage[all]{xy}
\usepackage{xr}
\DeclareMathAlphabet{\mathpzc}{OT1}{pzc}{m}{it}

\setmargins{32mm}{20mm}{14.6cm}{22cm}{1cm}{1cm}{1cm}{1cm}

\newtheorem{theorem}{Theorem}[section]

\newtheorem{corollary}[theorem]{Corollary}

\newtheorem{lemma}[theorem]{Lemma}
\newtheorem{theorem*}{Theorem}
\newtheorem{proposition*}[theorem*]{Proposition}
\newtheorem{corollary*}[theorem*]{Corollary}
\newtheorem{lemma*}[theorem*]{Lemma}

\theoremstyle{definition}
\newtheorem{definition}[theorem]{Definition}

\newtheorem{definition*}[theorem*]{Definition}

\theoremstyle{remark}

\newtheorem{remark}[theorem]{Remark}
\newtheorem{remarks}[theorem]{Remarks}

\newtheorem{example*}[theorem*]{Example}
\newtheorem{examples*}[theorem*]{Examples}
\newtheorem{remark*}[theorem*]{Remark}
\newtheorem{remarks*}[theorem*]{Remarks}
\newtheorem{exercise*}[theorem*]{Exercise}

\newcommand\id{\mathrm{id}}

\newcommand\ten{\otimes}

\newcommand\CC{\mathrm{C}}

\renewcommand\H{\mathrm{H}}

\newcommand\Z{\mathbb{Z}}
\newcommand\Q{\mathbb{Q}}

\newcommand\bB{\mathbb{B}}

\newcommand\bD{\mathbb{D}}

\newcommand\bF{\mathbb{F}}

\newcommand\bP{\mathbb{P}}

\newcommand\cL{\mathcal{L}}

\newcommand\cN{\mathcal{N}}

\newcommand\cU{\mathcal{U}}

\newcommand\cW{\mathcal{W}}
\newcommand\cX{\mathcal{X}}

\renewcommand\O{\mathscr{O}}

\newcommand\sC{\mathscr{C}}

\newcommand\sF{\mathscr{F}}

\newcommand\Def{\mathfrak{Def}}

\newcommand\fD{\mathfrak{D}}

\newcommand\fT{\mathfrak{T}}
\newcommand\fU{\mathfrak{U}}

\newcommand\fX{\mathfrak{X}}

\newcommand\DP{^{\mathrm{DP}}}

\newcommand\g{\mathfrak{g}}

\newcommand\fh{\mathfrak{h}}

\newcommand\fp{\mathfrak{p}}

\newcommand\Hom{\mathrm{Hom}}

\newcommand\Gal{\mathrm{Gal}}

\newcommand\im{\mathrm{Im\,}}

\newcommand\Spec{\mathrm{Spec}\,}

\newcommand\Set{\mathrm{Set}}

\newcommand\Cat{\mathrm{Cat}}

\newcommand\Grpd{\mathrm{Grpd}}

\newcommand\Lim{\varprojlim}

\newcommand\into{\hookrightarrow}
\newcommand\onto{\twoheadrightarrow}

\newcommand\xra{\xrightarrow}
\newcommand\xla{\xleftarrow}

\newcommand\bt{\bullet}
\newcommand\by{\times}

\newcommand\ddef{\mathrm{Def}}

\newcommand\et{\acute{\mathrm{e}}\mathrm{t}}

\newcommand\half{\frac{1}{2}}

\newcommand\cris{\rm{cris}}

\newcommand\Bcr{B_{\cris}}

\newcommand\Gr{\mathrm{Gr}}

\newcommand\ab{\mathrm{ab}}
\newcommand\cts{\mathrm{cts}}
\newcommand\gp{\mathrm{Gp}}

\newcommand\Zar{\mathrm{Zar}}
\newcommand\Hdg{\mathrm{Hdg}}
\newcommand\Fil{\mathrm{Fil}}

\begin{document}
\title{Galois actions on the  pro-$l$-unipotent fundamental group\thanks{The author was supported during this research by Trinity College, Cambridge and by the Isle of Man Department of Education.}}
\author{J.P.Pridham}
\maketitle

%
\tableofcontents

\section*{Introduction}
\addcontentsline{toc}{section}{Introduction}

In \cite{pi1},  the pro-$l$-unipotent \'etale fundamental group $\pi_1(X,\bar{x})\ten\Q_l$ of a variety $X$ was defined, and  the Weil conjectures used to prove restrictions on its structure when $X$ is smooth and defined over a finite field. The purpose of this paper is to show that this approach can also be used to describe the Galois action on  $\pi_1(X,\bar{x})\ten\Q_l$, with $X$ now defined over a number field $K$.

Let $G=\Gal(\bar{K}/K)$, and $\fp$ a prime of $\bar{K}$.  The main idea in Section \ref{one} is to define, for each $\fp$, a $G_{\fp}$-invariant decomposition on  $\cL(\pi_1(X,\bar{x}),\Q_l)$, the Lie algebra underlying the pro-$l$-unipotent fundamental group. This can be done for every prime $\fp \nmid l$ at which $X$ is of potentially good reduction. Weight restrictions on cohomology then give restrictions on the possible weight decompositions on this Lie algebra. In particular, when $X$ is proper, the $G_{\fp}$-action on $\pi_1(X,\bar{x})\ten\Q_l$  is completely determined by the $G_{\fp}$-action on cohomology. The case when $\fp \mid l$ is more difficult. A suitable weight decomposition can often be defined, but only after tensoring with $B_{\cris}^{\phi}$. This requires that $\cL(\pi_1(X_{\bar{K}},\bar{x}),\Q_l)$ be a potentially crystalline Galois representation, and the remainder of the paper is concerned with proving that  this is the case.

In Section \ref{pi1cris}, the pro-unipotent crystalline fundamental group is defined, together with a cosimplicial algebra which can be thought of as the crystalline homotopy type. The cohomology groups of the crystalline homotopy type are the crystalline cohomology groups of the scheme, and the pro-unipotent fundamental group can be recovered from the algebra.

The purpose of Section \ref{phodge} is to establish a comparison theorem between the crystalline homotopy type and the pro-$p$ \'etale homotopy type defined in \cite{pi1}. This just amounts to showing that the chain of isomorphisms of cohomology groups in \cite{Hop} can be replaced by a chain of quasi-isomorphisms of cosimplicial algebras. As a consequence of this comparison, we can compare fundamental groups (Corollary \ref{compthm}) when $X$ is of good reduction:
$$
\cL(\pi_1^{\cris}(X_k,x_s),K_0)\ten_{K_0}B \cong \cL(\pi_1^{\et}(X_{\bar{K}},\bar{x}_{\eta}),\Q_p)\ten_{\Q_p} B.
$$
so that $\cL(\pi_1^{\et}(X_{\bar{K}},\bar{x}_{\eta}),\Q_p)$ is indeed a crystalline Galois representation.

\section{The pro-$l$-unipotent fundamental group}\label{one}

As in \cite{pi1} Section \ref{pi1-nilplie}, define $\cN_{\Q_l}$ to be the category of finite-dimensional nilpotent Lie algebras over $\Q_l$, and $\widehat{\cN_{\Q_l}}$ to be the category of finitely generated (i.e. \mbox{$\dim L/[L,L] < \infty$}) pro-nilpotent Lie algebras over $\Q_l$.

Recall (ibid. Section \ref{pi1-nrep}) that, for a pro-finite group $\Gamma$, we define the $\Q_l$ Malcev Lie algebra $\cL(\Gamma,\Q_l) \in \widehat{\cN_{\Q_l}}$  to pro-represent the functor on $\cN_{\Q_l}$ given by  
$$
\g \mapsto \Hom_{\cts}(\Gamma, \exp(\g)).
$$
 Define the pro-$\Q_l$-unipotent completion  of $\Gamma$ by
$$
\Gamma \ten_{\hat{\Z}} \Q_l:=\exp(\cL(\Gamma,\Q_l)).
$$

Given a ground field $k$, a connected scheme $X_0/k$ and a  point $x \in X_0$,  let \mbox{$X=X_0\ten_k \bar{k}$,} and let $\bar{x} \in X$ be the point above $x$.   We wish to study the action of $\Gal(\bar{k}/k)$ on the group
$$
\pi_1(X,\bar{x})\ten\Q_l,
$$
for $l$ different from the characteristic of $k$. From ibid. Lemma \ref{pi1-mittag}, this fundamental group has tangent space $\H^1(X,\Q_l)$ and obstruction space $\H^2(X,\Q_l)$.

\begin{definition} Given a smooth proper variety $\bar{X}$ over $\bar{k}$, with $D$ a divisor locally of normal crossings, let \mbox{$X=\bar{X}-D$}, and fix a geometric point $\bar{x} \to X$. Define the ideal $K$ as the kernel of the map
$$
\cL(\pi_1(X,\bar{x}),\Q_l) \to \cL(\pi_1(\bar{X},\bar{x}),\Q_l).
$$
\end{definition}

\begin{lemma}\label{gdwgt}
The definition of $K$ is independent of the choice $\bar{X}$ of compactification of $X$.
\begin{proof}
Consider the Lie algebra $\cL(\pi_1(X,\bar{x}),\Q_l)/K$. This has tangent space $\im(\H^1(\bar{X},\Q_l) \to \H^1(X,\Q_l))$, and obstruction space $\H^2(X,\Q_l)$. Now, as in \cite{Hodge2} 3.2.4.1, and similarly to \cite{Mi} Corollary VI 5.3, we have the Gysin exact sequence
$$
\H^1(\bar{X},\Q_l) \to \H^1(X,\Q_l)) \to \H^0(\tilde{D},\Q_l(-1)),
$$
where $\tilde{D}$ is the normalisation of $D$. If we adopt the same conventions for weights as in \cite{Poids} (i.e. Frobenius in finite characteristic, and Hodge theory in characteristic zero) then we see that
$$
\im(\H^1(\bar{X},\Q_l) \to \H^1(X,\Q_l))= W_1\H^1(X,\Q_l),
$$
the component of weight one. Thus for any two choices of $\bar{X}$, our Lie algebra has the same tangent and obstruction spaces. 

As in \cite{Hodge2} 3.2.11 C, we my compare any two compactifications by taking morphisms to a third, and since these morphisms induce isomorphisms of the tangent and obstruction spaces of the corresponding Lie algebras, they must induce isomorphism on the Lie algebras themselves. 
\end{proof}
\end{lemma}

\begin{definition}\label{weight} In the situation above, define the increasing weight filtration on \mbox{$\cL:=\cL(\pi_1(X,\bar{x}),\Q_l)$} to be the filtration generated by the following conditions:
$$
W_{-1}\cL=\cL,\quad K\subset W_{-2}\cL, \quad [W_m\cL,W_n\cL] \subset W_{m+n}\cL.
$$
\end{definition}

Most of this paper is concerned with finding splittings of this weight filtration, i.e. weight decompositions. Note that if $X$ is smooth and proper then the weight filtration is just the commutator filtration $W_{-n}=\Gamma_n$.

\begin{remark} There are times when $X_0/k$ has no closed points, yet we may still define an \emph{outer} action of $\Gal(\bar{k}/k)$ on $\pi_1(X)$. We may regard $\pi_1(X)\ten \Q_l$ as being the hull of the deformation functor $\mathfrak{PB}(X)$ defined in \cite{pi1} Section \ref{pi1-nrep}. Alternatively, we may use tangential basepoints, as dealt with in ibid. Remark \ref{pi1-tangential}.
\end{remark}

\subsection{Finite characteristic}

Let $k=\bF_q$ be a finite field, with $X$  variety over $\bar{k}$,  either proper or smooth. There is a natural Frobenius action $F$ on $\pi_1(X,\bar{x})$, and hence on
$$
\cL:=\cL(\pi_1(X,\bar{x}),\Q_l).
$$

Let $\cL_n=\cL/(\Gamma_{n+1}\cL)$, recalling that $\Gamma_n$ is the $n$th commutator (i.e. $\Gamma_1\cL=\cL$, and $\Gamma_{n+1}\cL=[\cL,\Gamma_n\cL]$). Since $\cL_n$ is finite dimensional, we may consider Jordan decompositions (over $\bar{\Q_l}$) given by $F$.

\begin{lemma} 
On $\cL_n$, all of the eigenvalues of $F$ are algebraic numbers, all of whose conjugates have norm $q^{w/2}$, for some negative integer $w$ (dependent on $\alpha$).
\begin{proof}
We prove this by induction on $n$. For $n=1$, $\cL_1 \cong \H^1(X,\Q_l)^{\vee}$, so this statement is just \cite{Weil2} Corollaries (3.3.4)--(3.3.6).

Now, note that if $v$ is in the generalised $\lambda$-eigenspace for $F$, and $w$ in the generalised $\mu$-eigenspace, then
$[v,w]$ is in the generalised $\lambda\mu$-eigenspace, since
$$
(F-\lambda\mu)[v,w]=[(F-\lambda)v,Fw]+[\lambda v, (F-\mu)w],
$$
so $(F-\lambda\mu)^N[v,w]=0$ for $N \gg 0$.

Hence the statement of the Lemma holds for $\Gamma_n(\cL)/\Gamma_{n+1}(\cL)$ (by multiplicativity), since it holds for $\cL_1$.

Assume the statement holds for $\cL_{n-1}$. Consider the exact sequence
$$
0 \to \Gamma_n(\cL)/\Gamma_{n+1}(\cL) \to \cL_n \to \cL_{n-1} \to 0.
$$
Since the statement holds for the terms on either side, it holds for $\cL_n$.
\end{proof}
\end{lemma}

\begin{corollary} There is a weight decomposition $\cL_n=\bigoplus_i \cW_i(\cL_n)$ on the inverse system $\cL_n$. 
\begin{proof} We just decompose by the weights $i$ of the eigenvalues of norm $q^{i/2}$.
\end{proof}
\end{corollary}

\begin{theorem}\label{finitewgt} If $X$ is smooth, satisfying the conditions of Definition \ref{weight}, then this weight decomposition is a splitting of the weight filtration $W$.
\begin{proof}
It suffices to show that $\Gr^W_i\cL$ is pure of weight $i$ (in the sense of Frobenius eigenvalues). Indeed, it is enough for $\Gr_{\Gamma}^r\Gr^W_i\cL$ to be pure of weight $i$. But now observe that $\Gr_{\Gamma}^r\Gr^W_i\cL$ is spanned by commutators consisting of $i-r$ elements of $K/K\cap \Gamma_2$ and $2r-i$ elements of $\cL/(\Gamma_2+K)$.

Now, $\cL/(\Gamma_2+K)$ is dual to $W_1\H^1(X,\Q_l)$, which is pure of weight $1$, by Lemma \ref{gdwgt}, so $\cL/(\Gamma_2+K)$ is pure of weight $-1$. Meanwhile, $K/K\cap \Gamma_2$ is dual to $\H^1(X,\Q_l)/W_1\H^1(X,\Q_l)$, which is pure of weight $2$, so $K/K\cap \Gamma_2$ is pure of weight $-2$.

Therefore $\Gr_{\Gamma}^r\Gr^W_i\cL$ is pure of weight $i$, as required.
\end{proof}
\end{theorem}

\begin{corollary} For $X$ as above, there is a Frobenius equivariant isomorphism of inverse systems $\cL_n(\pi_1(X,\bar{x}),\Q_l) \cong \Gr^W\cL_n(\pi_1(X,\bar{x}),\Q_l)$.
\end{corollary}

\begin{theorem}\label{finitestr}
\begin{enumerate}
\item If  $\H^1(X,\Q_l)$ is pure of  weight $a>0$, then there is a Frobenius equivariant map
$$
L(\H^1(X,\Q_l)^{\vee}) \onto \cL(\pi_1(X,\bar{x}),\Q_l),
$$
inducing the identity on tangent spaces $\H^1(X,\Q_l)$.

\item If moreover $\H^2(X,\Q_l)$ is mixed of weights $\le 2a$, then there is a Frobenius equivariant isomorphism
$$
L(\H^1(X,\Q_l)^{\vee})/(\check{\cup}\,\H^2(X,\Q_l)^{\vee})) \cong \cL(\pi_1(X,\bar{x}),\Q_l),
$$
where $\check{\cup}$ is dual to the cup product.
\end{enumerate}
\begin{proof}
\begin{enumerate}
\item Since $\H^1(X,\Q_l)$ is pure of weight $a$, it follows that $\Gamma_n\cL=W_{-na}\cL$. It then follows that the composition
$$
\cW_{-a} \into \cL \onto \cL/\Gamma_2(\cL) \cong \H^1(X,\Q_l)^{\vee}
$$
is a Frobenius equivariant isomorphism, so we may take its inverse to give a map $\H^1(X,\Q_l)^{\vee} \to \cL$, extending to give
$$
L(\H^1(X,\Q_l)^{\vee}) \onto \cL(\pi_1(X,\bar{x}),\Q_l).
$$
 
\item This is essentially the same as the proof of \cite{pi1} Theorem \ref{pi1-nfrobhull}. Let $J$ be the kernel of the map above, and write $\H_1$ for $\H^1(X,\Q_l)^{\vee}$. It follows that $(J/[\H_1,J])^{\vee}$ is a universal obstruction space for $\cL$, which gives us a unique morphism
$$
(J/[\H_1,J])^{\vee} \into \H^2(X,\Q_l)
$$
of obstruction theories. These obstruction theories are both Frobenius equivariant, so this map must be. Now note that $J \subset \Gamma_2\, L(\H_1)$, so $J$ is of weights $\le -2a$, so $(J/[\H_1,J])^{\vee}$ is of weights $\ge 2a$. By hypothesis, $\H^2(X,\Q_l)$ is of weights $\le 2a$, so the inclusion above implies that $(J/[\H_1,J])^{\vee}$ must be pure of weight $2a$. 

Now, since  $[\H_1,J]$ is of weights $\le 3a$, the composition
$$
\cW_{-2a}\cap J \into J \onto J/[\H_1,J]
$$
must be an isomorphism, so $\cW_{-2a}\cap J$ generates the ideal $J$, and since $\H^1(X,\Q_l)$ is pure of weight $a$,
$$
\cW_{-2a} = {\bigwedge}^2 \H^1(X,\Q_l)^{\vee}.
$$
The inclusion of obstruction spaces therefore corresponds to  a Frobenius equivariant  map 
$$
\H^2(X,\Q_l)^{\vee} \onto \cW_{-2a}\cap J = {\bigwedge}^2 \H^1(X,\Q_l)^{\vee}\cap J,
$$
which, by comparing primary obstruction maps, must be dual to the cup product. Hence
$$
J= (\check{\cup}\,\H^2(X,\Q_l)^{\vee}),
$$
so
$$
L(\H^1(X,\Q_l)^{\vee})/(\check{\cup}\,\H^2(X,\Q_l)^{\vee})) \cong \cL(\pi_1(X,\bar{x}),\Q_l).
$$
\end{enumerate}
\end{proof}
\end{theorem}
 
\begin{remarks}
\begin{enumerate}
\item If $X$ is smooth and proper, then $\H^1(X,\Q_l)$ will be pure of weight $1$, and $\H^2(X,\Q_l)$ pure of  weight $2$, so both conditions will be satisfied.
\item If $X=\bar{X}-D$, for $\bar{X}$ smooth, proper and simply connected, and $D$ a divisor locally of normal crossings, then $\H^1(X,\Q_l)$ will be pure of weight $2$. Since $\H^2(X,\Q_l)$ is automatically of weights $\le 4$, both conditions of the theorem are satisfied.
\item If $X=\bar{X}-D$, with the Gysin map $\H^0(\tilde{D},\Q_l(-1)) \to \H^2(X,\Q_l)$ injective, then $\H^1(X,\Q_l)$ will be pure of weight $1$, so the first condition is satisfied.
\end{enumerate}
\end{remarks}

\subsection{Good reduction, $l \ne p$}\label{good}

In characteristic zero, the easiest case to deal with is when working over a local field. Let $V$ be a complete discrete valuation ring, with fraction field $K$, and finite residue field $k$ of characteristic $p$. Denote by $\bar{K}$ (resp. $\bar{k}$) the algebraic closure of $K$ (resp. $k$), and $\bar{V}$ the integral closure of $V$ in $\bar{K}$.  Write $G=\Gal(\bar{K}/K)$. $l$ will be a prime other than $p$.

Take a smooth proper scheme $\bar{X}_V$ over $V$, with $D_V \subset \bar{X}_V$ a divisor of simple normal crossings, and let $X_V=\bar{X}_V-D_V$. Let $x$ be a $V$-valued point of $X_V$, with corresponding points $x_{\eta} \in X_V(K),\, x_s \in X_V(k)$, and the geometric points $\bar{x}_{\eta} \in X_V(\bar{K}),\, \bar{x}_s \in  X_V(\bar{k})$. 

We begin with a smooth specialisation theorem for pro-$l$-unipotent fundamental groups.

\begin{theorem}
There is a  $G$-equivariant isomorphism
$$
\pi_1(X_{\bar{K}},\bar{x}_{\eta})\ten\Q_l \cong \pi_1(X_{\bar{k}},\bar{x}_s)\ten \Q_l,
$$
so that, in particular, the $G$ action is unramified, factoring through
$$
G \to \Gal(\bar{k}/k).
$$
\begin{proof}
Consider the maps
\begin{eqnarray*}
\pi_1(X_{\bar{K}},\bar{x}_{\eta})\ten \Q_l &\to& \pi_1(X_{\bar{V}},\bar{x}_{\eta})\ten \Q_l\\
\pi_1(X_{\bar{k}},\bar{x}_s)\ten \Q_l &\to& \pi_1(X_{\bar{V}},\bar{x}_s)\ten \Q_l.
\end{eqnarray*}
Comparing tangent and obstruction spaces over $\cN_{\Q_l}$, we have:
\begin{eqnarray*}
\H^i(X_{\bar{K}},\Q_l) &\leftarrow& \H^i(X_{\bar{V}},\Q_l)\\
\H^i(X_{\bar{k}},\Q_l) &\leftarrow& \H^i(X_{\bar{V}},\Q_l),
\end{eqnarray*}
for $i=1,2$. Since $\bar{V}$ is  strictly Henselian,  we may apply the smooth specialisation theorem for geometric fibrations (\cite{friedlanderfib} Propositions 2.3 and 2.5) to see that these are isomorphisms. Since \'etale maps of pro-representables on $\cN_{\Q_l}$ must be isomorphisms, it follows that we have isomorphisms on pro-$l$-unipotent fundamental groups. Note that these isomorphisms are $G$-equivariant.

It only remains to observe that the embedding $\bar{V} \to \bar{K}$ provides us with a canonical ($G$-equivariant) isomorphism
$$
\pi_1(X_{\bar{V}},\bar{x}_s) \to \pi_1(X_{\bar{V}},\bar{x}_{\eta}).
$$
\end{proof}
\end{theorem}

\begin{corollary} There is a $G$-invariant weight decomposition on $\cL_n(\pi_1(X_{\bar{K}},\bar{x}_{\eta}), \Q_l)$, splitting the weight filtration $W$. Moreover,
\begin{enumerate}
\item If  $\H^1(X_{\bar{K}},\Q_l)$ is pure of  weight $a>0$, then there is a Frobenius equivariant map
$$
L(\H^1(X_{\bar{K}},\Q_l)^{\vee}) \onto \cL(\pi_1(X_{\bar{K}},\bar{x}_{\eta}),\Q_l),
$$
inducing the identity on tangent spaces $\H^1(X_{\bar{K}},\Q_l)$.

\item If in addition $\H^2(X_{\bar{K}},\Q_l)$ is mixed of weights $\le 2a$, then there is a Frobenius equivariant isomorphism
$$
L(\H^1(X_{\bar{K}},\Q_l)^{\vee})/(\check{\cup}\,\H^2(X_{\bar{K}},\Q_l)^{\vee})) \cong \cL(\pi_1(X_{\bar{K}},\bar{x}_{\eta}),\Q_l),
$$
where $\check{\cup}$ is dual to the cup product.
\end{enumerate}
\begin{proof}
This follows from  Theorems \ref{finitewgt} and \ref{finitestr}.
\end{proof}
\end{corollary}

\begin{remark}
As in \cite{pi1} Remark \ref{pi1-htpy}, we see that, if $X$ is smooth and proper, then the  $\Q_l$-homotopy type is formal and the quasi-isomorphism is Galois equivariant.
\end{remark}  

\subsection{Potentially good reduction, $l \ne p$}\label{pgood}

Let $K$ be as before, and  $K'/K$ a finite Galois extension, with ring of integers $V'$. Fix a smooth scheme $X/K$,   with 
$$
X=\bar{X} -D,
$$
such that $\bar{X}$ is smooth and proper, and  $D$ a divisor of simple normal crossings, and fix a point $x \in X(K)$.

\begin{theorem}\label{pgdwgt}
 If there exists some $Y/V'$  smooth and proper, with $Z \subset Y$ a  divisor of simple normal crossings, such that $Y_{K'} =\bar{X}_{K'}$ and $Z_{K'} =D_{K'}$, then  there is a $G$-invariant weight decomposition on $\cL_n(\pi_1(X_{\bar{K}},\bar{x}), \Q_l)$, splitting the weight filtration $W$. Moreover,
\begin{enumerate}
\item If  $\H^1(X_{\bar{K}},\Q_l)$ is pure of  weight $a>0$, then there is a $G$-equivariant map
$$
L(\H^1(X_{\bar{K}},\Q_l)^{\vee}) \onto \cL(\pi_1(X_{\bar{K}},\bar{x}),\Q_l),
$$
inducing the identity on tangent spaces $\H^1(X_{\bar{K}},\Q_l)$.

\item If in addition $\H^2(X_{\bar{K}},\Q_l)$ is  of weight $\le 2a$, then there is a $G$-equivariant isomorphism
$$
L(\H^1(X_{\bar{K}},\Q_l)^{\vee})/(\check{\cup}\,\H^2(X_{\bar{K}},\Q_l)^{\vee})) \cong \cL(\pi_1(X_{\bar{K}},\bar{x}),\Q_l),
$$
where $\check{\cup}$ is dual to the cup product.
\end{enumerate}
\begin{proof}
Observe that if we replace $G$ by $G'=\Gal(\bar{K}/K')$, then  this is just Theorem \ref{gdwgt}. Take the $G'$-invariant weight decomposition $\cW_i$ from that theorem. It will suffice to show that this decomposition is in fact $G$-invariant.

Fix $g \in G$, and  note that $g\cW_*$ is the decomposition due to the Frobenius operator $g\tilde{F}g^{-1}$. Since $\cW_i$ is independent of the choice of Frobenius, if follows that $g\cW_i=\cW_i$. The remainder of the proof is now identical to Theorem \ref{gdwgt}.
\end{proof}
\end{theorem} 

\subsection{Potentially good reduction, $l=p$}\label{pgdl=p}

Having established  the structure of $\pi_1(X_{\bar{K}},\bar{x}_{\eta})\ten \Q_l$, it is natural to wish to understand  $\pi_1(X_{\bar{K}},\bar{x}_{\eta})\ten \Q_p$. By passing to a global field, and then localising at a different prime, we may describe the structure, but not the Galois action. If we could describe the Galois action on $\pi_1(X_{\bar{K}},\bar{x}_{\eta})\ten \Q_p$, then we would also be able to describe the action on its tangent space, $\H^1(X_{\bar{K}},\Q_p)$, but this requires $p$-adic Hodge theory.

Let $K$ be as before, and  $K'/K$ a finite Galois extension, with residue field $k'=\bF_q$, of characteristic $p$. Let $W'$ be the ring of Witt vectors $W=W(k')$, with   fraction field $K'_0$. Let $G=\Gal(\bar{K}/K)$, and $G'=\Gal(\bar{K}/K')$.   Fix a smooth scheme $X/K$,   with 
$$
X=\bar{X} -D,
$$
such that $\bar{X}$ is smooth and proper and  $D$ a divisor locally of normal crossings, and fix a point $x \in X(K)$. Assume that $\bar{X}$ is of good reduction over $K'$.

We wish to study the Galois structure of $\cL:=\cL(\pi_1(X_{\bar{K}},\bar{x}),\Q_p)$. We can no longer split the weight filtration here, since there is no appropriate Frobenius. Instead we have to use $p$-adic periods.

Consider the $K'_0$-Lie algebra
$$
D'_{\cris}(\cL):= \Lim D'_{\cris}(\cL_n)=\Lim (\cL_n \ten_{\Q_p} B_{\cris})^{G'},
$$
where $B_{\cris}$ is defined as in \cite{Tsuji}. 

This will automatically have a $\sigma$-linear Frobenius operator $\phi$, where $\sigma$ is the Frobenius operator on $K'_0$. It also has a  Hodge filtration $\Fil^*$ on $D'_{\cris}(\cL)\ten_{K'_0}K'$, and  an action of $G$, factoring through the finite group $G/G'$.  These structures all respect the Lie bracket,  while $\phi$  preserves the Hodge filtration.

The weight filtration given in Definition \ref{weight} extends naturally to give a weight filtration on $D'_{\cris}(\cL)$. Our aim is to show that this weight filtration can be split, using a weight decomposition as defined in Section \ref{gdwgt}.

Assume that $\cL$ is an inverse system of crystalline $G'$-modules.  Section \ref{phodge} gives conditions for this hypothesis to hold. Having made this assumption, we are able to  proceed as before.

Now, observe that 
$$
\Gr^W_{i} \Gr_{\Gamma}^1 D'_{\cris}(\cL) = \Gr^W_{i} D'_{\cris}((\H^1(X_{\bar{K}},\Q_p))^{\vee}),
$$
and recall that
\begin{eqnarray*}
\H^1(\bar{X}_{\bar{K}},\Q_p) &\onto&  \Gr^W_1 \H^1(X_{\bar{K}},\Q_p),\\
\Gr^W_2 \H^1(X_{\bar{K}},\Q_p) &\into& \H^0(\tilde{D}_{\bar{K}},\Q_p(-1)),
\end{eqnarray*}
with all other graded weights being zero.

It follows automatically that $D'_{\cris}(\H^0(\tilde{D}_{\bar{K}},\Q_p(-1)))$ is a direct sum of  vector spaces $K'_0[-1]$, on which $\phi=p\sigma$,  and $G$ acts as on $K'_0$,  the Hodge filtration being concentrated in degree one.
 From the crystalline conjecture (proved in \cite{Hop} or \cite{Tsuji}), it also follows that 
$$
D'_{cr}(\H^1(\bar{X}_{\bar{K}},\Q_p)) \cong \H^1_{\cris}(\bar{X}_{k'}/W')\ten K'_0.
$$

Next, consider the action of Frobenius on $D'_{\cris}(\cL)$. Although $\phi$ is only $\sigma$-linear, we will have $\phi^r$ linear, where $q =p^r$. We may therefore consider the eigenvalues of Frobenius. It is clear that this action is pure of weight $-2$ on 
$
\Gr^W_{-2} \Gr_{\Gamma}^1 D'_{\cris}(\cL).
$

From \cite{KM}, for any smooth projective curve $Y/k'$, the Frobenius action on
$$
\H^1_{\cris}(Y/W')\ten K'_0
$$
is pure of weight $1$. Therefore, for any smooth projective curve $Y/K'$ of good reduction, we see that the Frobenius action on
$$
D_{\cris}(\H^1_{\et}(Y_{\bar{K}},\Q_p))
$$
is pure of weight $1$.

By considering the Albanese variety, we see that this must hold for any smooth proper variety $Y/K'$ of good reduction, so in particular holds for $\bar{X}$. Now, 
$$
\Gr^W_{-1} \Gr_{\Gamma}^1 D'_{\cris}(\cL) \into D'_{\cris}((\H^1_{\et}(\bar{X}_{\bar{K}},\Q_p))^{\vee}),
$$
so the Frobenius action on 
$$
\Gr^W_{-1} \Gr_{\Gamma}^1 D'_{\cris}(\cL)
$$
is pure of weight $-1$. 

\begin{lemma}
There is a $(G,\phi)$-invariant weight decomposition $\cW_*$  on $D'_{\cris}(\cL_n)$, splitting the weight filtration $W$. 
\begin{proof}
The  proof follows almost exactly as for Theorem \ref{gdwgt}, one minor difference being that the Frobenius decomposition is defined by $\phi^r$, not $\phi$. Since $\phi$ commutes with $\phi^r$, it preserves the weight decomposition. The action of $G$ commutes with that of  $\phi^r$ (though not $\phi$), so preserves $\cW$.
\end{proof}
\end{lemma}

\begin{remark}
In general,  $\cW_*$ is incompatible with the Hodge filtration $\Fil^*$, as the latter is not preserved by $\phi$.
\end{remark}

\begin{theorem}\label{crpgdwgt}
If the pair $(\bar{X}_{K'}, D_{K'})$ is of good reduction for some finite Galois extension $K'/K$, and $\{\cL_n\}$ is a system of potentially crystalline Galois representations,  then there is a $G$-invariant weight decomposition on $\cL_n(\pi_1(X_{\bar{K}},\bar{x}), \Q_p)\ten_{\Q_p}B_{\cris}^{\phi}$, splitting the weight filtration $W$, where $B_{\cris}^{\phi}$ is the $\phi$-invariant subring of $B_{\cris}$. Moreover,
\begin{enumerate}
\item If  $D_{\cris}'\H^1(X_{\bar{K}},\Q_p)$ is pure of  weight $a>0$, then there is a canonical $G$-equivariant map
$$
L(\H^1(X_{\bar{K}},\Q_p)^{\vee})\ten_{\Q_p}B_{\cris}^{\phi} \onto \cL(\pi_1(X_{\bar{K}},\bar{x}),\Q_p)\ten_{\Q_p}B_{\cris}^{\phi},
$$
inducing the identity on tangent spaces $\H^1(X_{\bar{K}},\Q_p)\ten_{\Q_p}B_{\cris}^{\phi}$.

\item If in addition $D_{\cris}'\H^2(X_{\bar{K}},\Q_p)$ is  of weight $\le 2a$, then there is $G$-equivariant isomorphism
$$
B_{\cris}^{\phi}\ten_{\Q_p}L(\H^1(X_{\bar{K}},\Q_p)^{\vee})/(\check{\cup}\,\H^2(X_{\bar{K}},\Q_p)^{\vee})) \cong B_{\cris}^{\phi}\ten_{\Q_p}\cL(\pi_1(X_{\bar{K}},\bar{x}),\Q_p),
$$
where $\check{\cup}$ is dual to the cup product.
\end{enumerate}
\end{theorem}
\begin{proof}
This is very similar to Theorem \ref{gdwgt}. The weight decomposition is just $(\cW_*\ten_{K_0'} B_{\cris})^{\phi}$.

Letting $\cL:=\cL(\pi_1(X_{\bar{K}},\bar{x}),\Q_p)$ and $\cL^{\ab}:=\cL/[\cL,\cL]$, the first condition amounts to saying that  $D_{\cris}'(\cL^{\ab})$ is pure of weight $-a$. This gives a unique $\phi$-equivariant section of $D_{\cris}'(\cL) \to  D_{\cris}'(\cL^{\ab})$, and hence a canonical $(G,\phi)$-equivariant section of
$$
\cL \ten_{\Q_p}B_{\cris} \to \cL^{\ab} \ten_{\Q_p}B_{\cris},
$$
which proves the first part.

Now  write $V$ for $D_{\cris}'(\H^1(X_{\bar{k}},\Q_l))^{\vee}$, and let $J$ be the kernel of the map $L(V) \to D_{\cris}'(\cL)$ defined above.
Considering the extension
$$
e:(L(V)/[V,J]) \ten_{K_0'} B_{\cris} \onto \cL\ten_{\Q_p} B_{\cris}, 
$$
it follows that $(J/[V,J])^{\vee}\ten_{K_0'} B_{\cris}  \cong \H^2(\cL,\Q_l)\ten_{\Q_p} B_{\cris}$ is the universal obstruction space for $\cL\ten_{\Q_p} B_{\cris}$ over $B_{\cris}$, which gives us a unique $(G, \phi)$-equivariant morphism
$$
(J/[V,J])^{\vee}\ten_{K_0'} B_{\cris}  \into \H^2(X_{\bar{k}},\Q_l)\ten_{\Q_p} B_{\cris}
$$
of obstruction theories. Taking $G'$-invariants  gives the $\phi$-equivariant map
$$
(J/[V,J])^{\vee} \into D_{\cris}'(\H^2(X_{\bar{k}},\Q_l)).
$$
The remainder of the proof is as for Theorem \ref{finitestr}.
\end{proof}

\begin{remarks}
In Corollary \ref{crys} it will be shown that, when $D$ is a divisor of simple normal crossings, 
$\cL$ is a potentially crystalline representation. This then gives more information about the Galois action on  $\cL(\pi_1^{\et}(X_{\bar{K}},\bar{x}_{\eta}),\Q_p)$  than Olsson's  later results (\cite{olssonhodge} Theorem 1.13) obtained independently, which show that the Galois action on $\cL\ten\ten_{\Q_p}B_{\cris}$ is determined by the action on $\H^*(X,\Q_p)\ten_{\Q_p}B_{\cris}$   when $X$ is proper, without giving an explicit description.

In particular, note that all conditions of the theorem are satisfied by $X=\bP^1-\{0,1,\infty\}$, as studied in \cite{droite} and \cite{hainmat}, since its cohomology is pure. 
\end{remarks}

\subsection{Global fields}

We are now in a position to consider global fields. Let $K/\Q$ be a number field. Let
$$
X=\bar{X}-D,
$$
where $\bar{X}/K$ is a smooth proper variety, and $D$ is a divisor locally of normal crossings. Let $x \in X(k)$. For each prime $l$, we wish to consider the action of $G=\Gal(\bar{K}/K)$ on the pro-$l$-unipotent fundamental group
$$
\pi_1(X_{\bar{K}},\bar{x})\ten \Q_l.
$$
As in \cite{Poids}, we have a weight decomposition for each prime:

Now assume that the divisor $D \subset \bar{X}$ is of simple normal crossings.

\begin{theorem} For each prime $\fp \nmid l$ of $\bar{K}$ at which $(\bar{X},D)$ has potentially good reduction,  there is a weight decomposition ${}_{\fp}\cW_* (\cL)$ of $\cL:=\cL(\pi_1(X_{\bar{K}},\bar{x}), \Q_l)$, splitting the weight filtration of Definition \ref{weight}. These decompositions are conjugate under the action of $G$, i.e.
$$
g( {}_{\fp}\cW_*(\cL))= {}_{g\fp}\cW_*(\cL).
$$

For each prime $\fp \mid l$ at which $(\bar{X},D)$ has potentially good reduction and $\cL$ is a potentially crystalline $G_{\fp}$ representation, there is a  weight decomposition ${}_{\fp}\cW_* (\cL\ten_{\Q_p}B_{\cris}^{\phi})$ of $\cL\ten_{\Q_p}B_{\cris}^{\phi}$, splitting the weight filtration. These decompositions are conjugate.
\end{theorem}
\begin{proof}
We prove this for $\fp \nmid l$; the case $\fp \mid l$ is similar.  Localisation at $\fp$ gives  a $G_{\fp}$-equivariant isomorphism
$$
\pi_1(X_{\bar{K}},\bar{x})\ten \Q_l \xla{\theta_{\fp}}  \pi_1(X_{\bar{K}_{\fp}},\bar{x}_{\fp}) \ten \Q_l,
$$
preserving the weight filtration. Since $(\bar{X},D)$ is of  potentially good reduction at $\fp$, by Theorem \ref{pgdwgt} (or Theorem \ref{crpgdwgt} when $\fp \mid l$) there is a  $G_{\fp}$-invariant  splitting $\cW_i$ of the weight filtration on the  Lie algebra associated to the latter fundamental group. Set
$$
{}_{\fp}\cW_i(\cL):= \theta_{\fp}(\cW_i(\cL(\pi_1(X_{\bar{K}_{\fp}},\bar{x}_{\fp}), \Q_l))).
$$

To see that these are conjugate, consider the commutative diagram
$$
\begin{CD}
\pi_1(X_{\bar{K}},\bar{x})\ten \Q_l @<{\theta_{\fp}}<<  \pi_1(X_{\bar{K}_{\fp}},\bar{x}_{\fp}) \ten \Q_l\\
@VgVV							@VgVV\\
\pi_1(X_{\bar{K}},\bar{x})\ten \Q_l @<{\theta_{g\fp}}<<  \pi_1(X_{\bar{K}_{g\fp}},\bar{x}_{g\fp}) \ten \Q_l.
\end{CD}
$$
By uniqueness of the weight decomposition over local fields, we must have $g(\cW_i)=\cW_i$ on the right, so
$$
g({}_{\fp}\cW_*(\cL))= {}_{g\fp}\cW_*(\cL).
$$
\end{proof}

\section{The crystalline fundamental group}\label{pi1cris}

In this section, we set about defining a pro-unipotent crystalline fundamental group, analogous to the pro-unipotent \'etale fundamental group defined in \cite{pi1}.

\subsection{Smooth proper schemes}

Given a  perfect field $k$, let $W$ be the ring of Witt vectors $W=W(k)$, with DP-structure $\gamma$ and  fraction field $K$. We wish to define, for any  connected scheme $X/k$, and any point $x \in X(k)$, a pro-nilpotent crystalline fundamental group
$$
\pi_1^{\cris}(X,x)\ten K.
$$

\begin{definition} Given a  smooth group scheme $H/W$, let $H_n =H\ten_{W}W_n$. Define a constructible $H$-torsor $\bD$ on $X$ to be an inverse system of  principal $H$-torsors $\bD_n$  on the crystalline sites $X/W_n$,
such that 
$$
\bD_n=i_n^{-1} \bD_{n+1},
$$
for $i_n:(X/W_n)_{\cris} \to (X/W_{n+1})_{\cris}$,
so we may regard  $\bD$ as
$$
\bD=\Lim \bD_n.
$$

The principal $H$-torsor is then given by the sheaf
$$
(U \into T) \mapsto H(T),
$$
for $U \subset X$, and $U \to T$ a DP-thickening.

\end{definition}

\begin{definition}
Given a smooth group scheme $G/K$, define a constructible $G$-torsor $\bB$ to be 
$$
\bB=\bD\ten_W K,
$$
for some group scheme $H/W$ with $G=H\ten_{W}K$, and some constructible $H$-torsor $\bD$. To make this explicit, a constructible $G$-torsor is a pair
$
(H,\bD),
$
with $H/W$ as above, and  if we have a map $H_1 \to H_2$ of smooth group schemes over $W$, becoming the identity map on $G$ when tensored with $K$, then we regard the $H_2$-torsor
$$
T \mapsto \bD_1(T)\by^{H_1(T)}H_2(T)
$$
as being equivalent to the $H_1$-torsor $(H_1,\bD_1)$.
\end{definition}

We are now in a position to define a functor on $\cN_K$ which $\pi_1^{\cris}(X,x)\ten K$ should pro-represent. Given any $\g \in \cN_K$, observe that $\exp(\g)$ is represented by a smooth affine group scheme over $K$. Explicitly, we may construct this group scheme by considering the polynomial algebra over $K$ generated by $\g^{\vee}$, with comultiplication coming from the Campbell-Baker-Hausdorff formula, which will converge, since $\g$ is nilpotent.

\begin{definition}
Define the functor
$$
\fT(X): \cN_K \to \Grpd
$$
to associate to $\g$ the groupoid of constructible $\exp(\g)$-torsors on $X$. 

Define the groupoid 
$$
\fT_x(X)(\g)
$$
to be the fibre of
$$
\fT(X)(\g) \to \fT(x)(\g)
$$
over $(\exp(\g\ten\O_{X/W}),\id)$.
\end{definition}

The crucial difference between the groupoids $\fT_x(X)(\g)$ and $\fT(X)(\g)$ is that only those morphisms which give the identity on $\bB_x$ are permitted in the former.

\begin{definition}
Define
$$
\cL(\pi_1^{\cris}(X,x),K)
$$
to pro-represent the functor $\fT_x(X)$. Note that this functor is clearly homogeneous, and has trivial automorphisms, so will be pro-representable provided the relevant tangent space (which will be $\H^1_{\cris}(X/W)\ten K$) is finite dimensional. Define
$$
\pi_1^{\cris}(X,x)\ten K := \exp(\cL(\pi_1^{\cris}(X,x),K)).
$$
\end{definition}

\begin{remark}
Since exponentiation gives an equivalence between $\cN_K$ and the category of unipotent algebraic $K$-groups, Tannakian duality implies that linear representations of $\pi_1^{\cris}(X,x)\ten K$ correspond to nilpotent isocrystals on $X$. This means that our definition  of $\pi_1^{\cris}$ agrees with \cite{Shiho} Definition  4.1.5.
\end{remark}

\begin{definition}
Let 
$$
\mathfrak{FHS}(X):\cN_{\Q_l} \to \Cat
$$
associate to $\g$ the category of faithful constructible homogeneous $\exp(\g)$-sheaves on $X$. Observe that
$\mathfrak{T}(X)$ is the fibre of $\mathfrak{FHS}(X)$ over the trivial $\exp(0)$-sheaf $1$ (the constant sheaf of singleton sets).
\end{definition}

In order to find a suitable SDC to govern this problem, we will follow \cite{B} 5.28 and 5.26. Assume that we have a smooth formal scheme $\fX/W$, and a closed immersion $i:X \to \fX$.

 Take an affine cover $\fU_{\alpha}$ of $\fX$, let 
$$
U_{\alpha}=\fU_{\alpha}\ten k, \quad X'= \coprod_{\alpha}U_{\alpha} 
$$
 and, for each $n$, form the $W_n$-PD envelope 
$$
U_{\alpha} \into U_{\alpha}^{\mathrm{DP},n}
$$
associated to the embedding
$$
U_{\alpha} \into \fU_{\alpha}\ten W_n,
$$
and let
$$
{X'}^{\mathrm{DP},n}= \coprod_{\alpha}U_{\alpha}^{\mathrm{DP},n}.
$$
Similarly, form the $W_n$-PD envelope 
$$
X \into X^{\mathrm{DP},n},
$$
associated to the embedding $X \into \fX\ten W_n$.

%

Consider the maps (as in \cite{B} 5.26):
$$\label{criscompare}
\xymatrix{
X'_{\cris}|_{{X'}^{\mathrm{DP},n}} \ar[r]^{\phi} \ar[d]^j \ar[rrd]^(0.6){v} & {X'}^{\mathrm{DP},n}_{\Zar} \ar@{<->}[d] \ar[rrd]^(0.4){v}   \\
X'_{\cris} \ar[r]^{u_{X'/W_n}}\ar[rrd]^(0.6){v} & X'_{\Zar} \ar[rrd]^(0.4){v} &X_{\cris}|_{{X}^{\mathrm{DP},n}} \ar[r]_{\phi} \ar[d]_>>>{j}  & {X}^{\mathrm{DP},n}_{\Zar} \ar@{<->}[d]  \\
&&X_{\cris} \ar[r]_{u_{X/W_n}} & X_{\Zar},
}
$$

Now, adopting suitable $p$-adic conventions, we have a functorial comonadic adjunction
$$
\xymatrix@1{
\mathfrak{FHS}(X_{\cris}) \ar@<1ex>[r]^-{w^{-1}} & \mathfrak{FHS}({X'}_{\cris}|_{{X'}^{\mathrm{DP}}}) \ar@<1ex>[l]^-{w_*}_-{\bot},
}
$$
where $w=uj$.
The adjunction $j^{-1} \dashv j_{*}$ is comonadic, since ${X}_{\cris}|_{{X}^{\mathrm{DP}}} \xra{j} X_{\cris}$ is a covering (as in \cite{B} 5.28). The adjunction $v^{-1} \dashv v_{*}$ is also comonadic, since the maps $X' \to X$ are coverings.  Note that $w=vj=jv$, and that the maps $v^{-1},v_*,j^{-1},j_*$ all commute with each other (whenever the composition can be defined),  since $v^{-1}$ and $j^{-1}$  always commute with everything, and direct images commute with each other.

Next, observe that the maps
$$
\xymatrix@1{
\mathfrak{FHS}(X'_{\cris}|_{{X'}^{\mathrm{DP}}}) \ar@<1ex>[r]^-{\varphi_*} & \mathfrak{FHS}(X'_{\Zar}) \ar@<1ex>[l]^-{\varphi^*}
}
$$
give an equivalence of categories. The second  category above has uniformly trivial deformation theory,  as $\exp(\g)$ is smooth, and torsors on an affine scheme are trivial.

Hence deformations are described by the SDC
$$
S^n_X(\g)=\Hom_{\mathfrak{FHS}({X'}_{\cris}|_{{X'}\DP})(\g)}(w^{-1}\exp(\g\ten \O_{X/W}),(w^{-1}w_*)^n w^{-1}\exp(\g\ten \O_{X/W})).
$$

Now, write 
$$
\CC^n(X,\sF):= \Gamma({X'}_{\cris}|_{{X'}}\DP, (w^{-1}w_*)^n w^{-1}\sF),
$$
 for sheaves $\sF$ on $(X/W)_{\cris}$, and
$$
\sC^n(\sF):=(w_* w^{-1})^{n+1}\sF,
$$
so that $\CC^n(X,\sF)=\Gamma(X, \sC^n(\sF))$. \label{crishtpy}

We then  have  isomorphisms
\begin{eqnarray*}
S^n(\g) \cong \Gamma({X'}_{\cris}|_{{X'}\DP},(w^{-1}w_*)^n w^{-1}\exp(\g\ten \O_{X/W})) &\cong& \CC^n(X, \exp(\g\ten \O_{X/W})),\\
(b \mapsto g\cdot b) &\mapsfrom& g,
\end{eqnarray*}
the latter being a cosimplicial complex of groups, becoming an SDC via the Alexander-Whitney cup product. Finally, note that
$$
\CC^n(X, \exp(\g\ten \O_{X/W})) \cong \exp(\CC^n(X, \O_{X/W})\ten_W \g),
$$
so following \cite{pi1} Remark \ref{pi1-htpy}, we regard
$$
\CC^n(X, \O_{X/W})\ten_W K
$$
as being the crystalline homotopy type of $X$. 

\begin{remark}
This is essentially the same as the crystalline homotopy type defined in \cite{hainkim}, which is taken as the Thom-Whitney  algebra of the cosimplicial algebra $\CC^{\bt}_{\Zar}(X,\Omega^{\bt}_{X/W})\ten_WK$. If $D$ denotes the denormalisation functor from DG algebras to cosimplicial algebras, we have quasi-isomorphisms
$$
\CC^{\bt}(X, \O_{X/W})\ten_W K \to \CC^{\bt}(X,D\Omega^{\bt}_{X/W} )\ten_W K \leftarrow \CC^{\bt}_{\Zar}(X,D\Omega^{\bt}_{X/W})\ten_WK,
$$
giving the required equivalence. Since the homotopy categories of cosimplicial and DG algebras are equivalent, for convenience we work exclusively  with cosimplicial algebras.   
\end{remark}

\begin{lemma}
If $\H^i_{\cris}(X/W)$ and $\H^{i-1}_{\cris}(X/W)$ are finitely generated, then the cohomology groups of $S^{\bullet}_X$ are isomorphic to 
$$
\H^*_{\cris}(X/W)\ten_{W} K.
$$ 
\begin{proof}
As for \cite{pi1}  Lemma \ref{pi1-mittag}.
\end{proof}
\end{lemma}

\begin{theorem}\label{crispointed}
Assume that $x \in X(k)$ lifts to $\tilde{x} \in \fX(W)$. This enables us to define the morphism
$$
S^{\bullet}_X \xra{i_x^*}  S^{\bullet}_{x}
$$
of SDCs. 
Define the SDC
$$
S^n_{X,x} \subset S^n_X
$$
to be the fibre of $i_x^*$ over the unique $0$-valued point of $S^n_{x}$,
then  
$$
\Def_{S_{X,x}}
$$
is equivalent to the discrete groupoid  pro-represented by $\cL(\pi_1^{\cris}(X,x),K)$.

The cohomology of $S^{\bullet}_{X,x}$ is
$$
\H^i_{\cris}((X,x)/W)\ten K\cong \left\{\begin{matrix} 0 & i=0 \\ \H^i_{\cris}(X/W)\ten K & i \ne 0, \end{matrix} \right.
$$
provided that the $\H^i_{\cris}(X/W)$ are finitely generated.
\end{theorem}

Observe that, since these constructions were functorial, the absolute Frobenius endomorphism $F$ on $X$ gives rise to a Frobenius endomorphism on $S^{\bullet}_X$. If $k$ is the finite field $\bF_{q}$, with $q=p^r$, then $F^r$ is a morphism of $W$-schemes, and we have the  commutative diagram
$$
\begin{CD}
\Spec k @>{x}>> X \\
@VVV @VV{F^r}V \\
\Spec k @>{x}>> X.
\end{CD}
$$
 This then provides a Frobenius endomorphism of 
$$
S^n_{X,x},
$$
which will agree with the corresponding Frobenius endomorphism on 
$$
\pi_1^{\cris}(X,x)\ten K.
$$

Assume that all the eigenvalues of Frobenius $F^r$ acting on the cohomology group $\H^1_{\cris}(X/W)\ten K$  are algebraic numbers $\alpha$, and for each $\alpha$, there exists a weight $n$, such that all complex conjugates of $\alpha$ have norm $q^{n/2}$. This provides us with a weight decomposition 
$$
\H^1_{\cris}(X/W)\ten K=\bigoplus_n \cW_n \H^1_{\cris}(X/W)\ten K.
$$ 
In particular, this will be true if $X$ is smooth and projective (by \cite{KM}). This will enable the definition of a weight decomposition on $\cL(\pi_1^{\cris}(X,x), K)$.

\begin{theorem}\label{crisfrobhull}
Under the above hypotheses on eigenvalues of Frobenius, there is an isomorphism
$$
\cL(\pi_1^{\cris}(X,x),K) \cong  L(\H^1_{\cris}(X/W)\ten K^{\vee})/(f(\H^2_{\cris}(X/W)\ten K^{\vee})),
$$
where $L(V)$ is the free pro-nilpotent Lie algebra (as in \cite{pi1} \ref{pi1-nfreelie}), and  
$$
f:\H^2_{\cris}(X/W)\ten K^{\vee} \to \Gamma_2L(\H^1_{\cris}(X/W)\ten K^{\vee})
$$ 
preserves the (Frobenius) weight decompositions. The isomorphism also preserves the Frobenius decompositions of the finite quotients of  $\cL(\pi_1^{\cris}(X,x),K)$.

Moreover, the quotient map
$$
f:\H^2_{\cris}(X/W)\ten K^{\vee} \to \Gamma_2/\Gamma_3 \cong {\bigwedge}^2(\H^1_{\cris}(X/W)\ten K^{\vee})
$$
is dual to half the  cup product
$$
\H^1_{\cris}(X/W)\ten K \by \H^1_{\cris}(X/W)\ten K \xra{\half\cup} \H^2_{\cris}(X/W)\ten K.
$$
\begin{proof}
As for \cite{pi1} Theorem \ref{pi1-nfrobhull}.
\end{proof}
\end{theorem}

\begin{corollary}\label{crisquad}
If $X$ is smooth and projective, then 
$$
\cL(\pi_1^{\cris}(X,x),K)
$$
is quadratically presented. 

In fact,
$$
\cL(\pi_1^{\cris}(X,x),K) \cong  L(\H^1_{\cris}(X/W)\ten K^{\vee})/(\check{\cup}\,\H^2_{\cris}(X/W)\ten K^{\vee})),
$$
where $\check{\cup}$ is dual to the cup product. This isomorphism is, moreover, Frobenius equivariant.
\begin{proof}
This follows as for Theorem \ref{finitestr}, using the purity results of \cite{KM}.
\end{proof}
\end{corollary}

\begin{remark} Of course, in this case the Frobenius decomposition provides a decomposition on the weight filtration (which in this case is just the commutator filtration). The same result will hold for any smooth proper scheme whose  Frobenius action on $\H^1_{\cris}$ is pure of weight $i$.
\end{remark}

\subsection{Smooth non-proper schemes
}

There is a similar approach for the log-crystalline fundamental group. Adopting the conventions of \cite{Hop},  let $X^{o}=X - D$, for $X/k$ a smooth proper scheme, and $D$ a divisor locally of normal crossings. We wish to define, for $x \in X^o(k)$, a pro-nilpotent log-crystalline fundamental group
$$
\pi_1^{\cris}(X^o,x)\ten K.
$$

We define constructible torsors as before, replacing the crystalline site by the log-crystalline site throughout.

\begin{definition}
Define the functor
$$
\fT(X^o): \exp(\g) \to \Grpd
$$
to associate to $\g$ the groupoid of constructible $\exp(\g)$-torsors on the log-crystalline site of $X^o$.

Define the groupoid 
$$
\fT_x(X^o)(\g)
$$
to be the fibre of
$$
\fT(X^o)(\g) \to \fT(x)(\g)
$$
over $(\exp(\g),\id)$.
\end{definition}

\begin{definition}
Define
$$
\cL(\pi_1^{\cris}(X^o,x),K)
$$
to pro-represent the functor $\fT_x(X^o)$. 
\end{definition}

We now wish  to find a suitable SDC to govern this problem. Assume that we have a smooth formal scheme $\fX/W$, with a smooth formal divisor $\fD/W$ and  closed immersions $i:X \to \fX$, $D \to \fD$. 

 Take an affine cover $\fU_{\alpha}$ of $\fX$,  and, as before, form the  schemes 
$$
X^{\mathrm{DP},n},\quad \quad D^{\mathrm{DP},n},\quad \quad   U_{\alpha}^{\mathrm{DP},n} \quad \text{ and } \quad {X'}^{\mathrm{DP},n}.   
$$

In the notation of the previous section, We now  have the functorial comonadic adjunction
$$
\xymatrix@1{
\mathfrak{FHS}(X^o_{\cris}) \ar@<1ex>[r]^-{w^{-1}} & \mathfrak{FHS}({(X^o)'}_{\cris}|_{{(X^o)'}\DP}) \ar@<1ex>[l]^-{w_*}_-{\bot},
}
$$
the second  category  having uniformly trivial deformation theory.

Hence deformations are described by the SDC
$$
S^n_{X^o}(\g)=\Hom_{\mathfrak{FHS}({(X^o)'}_{\cris}|_{{(X^o)'}\DP})(\g)}(w^{-1}\exp(\g\ten \O_{X^o/W}),(w^{-1}w_*)^n w^{-1}\exp(\g\ten \O_{X^o/W})).
$$

Now, write 
$$
\CC^n(X^o,\sF):= \Gamma( {(X^o)'}_{\cris}|_{{(X^o)'}\DP}   , (w^{-1}w_*)^n w^{-1}\sF),
$$
 for sheaves $\sF$ on $(X^o/W)_{\cris}$, and
$$
\sC^n(\sF):=(w_* w^{-1})^{n+1}\sF,
$$
so that $\CC^n(X^o,\sF)=\Gamma(X^o, \sC^n(\sF))$. \label{logcrishtpy}

We then  have  isomorphisms
$$
S^n(\g) \cong \Gamma({(X^o)'}_{\cris}|_{{(X^o)'}\DP},(w^{-1}w_*)^n w^{-1}\exp(\g\ten \O_{X^o/W})),
$$
and
\begin{eqnarray*}
S^n(\g) &\cong& \CC^n(X^o, \exp(\g\ten \O_{X^o/W})),\\
(b \mapsto g\cdot b) &\mapsfrom& g,
\end{eqnarray*}
the latter being a cosimplicial complex of groups, becoming an SDC via the Alexander-Whitney cup product. Since
$$
\CC^n(X^o, \exp(\g\ten \O_{X^o/W})) \cong \exp(\CC^n(X^o, \O_{X^o/W})\ten_W \g),
$$
 we regard
$$
\CC^n(X^o, \O_{X/W})\ten_W K
$$
as being the log-crystalline homotopy type of $X^o$.

\begin{lemma}
If $\H^i_{\cris}(X^o/W)$ and $\H^{i-1}_{\cris}(X^o/W)$ are finitely generated, then the cohomology groups of $S^{\bullet}_{X^o}$ are isomorphic to 
$$
\H^*_{\cris}(X^o/W)\ten_{W} K.
$$ 
\begin{proof}
As for \cite{pi1} Lemma \ref{pi1-mittag}.
\end{proof}
\end{lemma}

\begin{theorem}\label{logcrispointed}
Assume that $x \in X^o(k)$ lifts to $\tilde{x} \in \fX(W)-\fD(W)$. This enables us to define the morphism
$$
S^{\bullet}_{X^o} \xra{i_x^*}  S^{\bullet}_{x}
$$
of SDCs. 
Define the SDC
$$
S^n_{X^o,x} \subset S^n_{X^o}
$$
to be the fibre of $i_x^*$ over the unique $0$-valued point of $S^n_{x}$,
then  
$$
\Def_{S_{X^o,x}}
$$
is equivalent to the discrete groupoid  pro-represented by $\cL(\pi_1^{\cris}(X^o,x),K)$.

The cohomology of $S^{\bullet}_{X^o,x}$ is
$$
\H^i_{\cris}((X^o,x)/W)\ten K\cong \left\{\begin{matrix} 0 & i=0 \\ \H^i_{\cris}(X^o/W)\ten K & i \ne 0, \end{matrix} \right.
$$
provided that  the $\H^i_{\cris}(X^o/W)$ are finitely generated.
\end{theorem}

Once again, if $k$ is finite, $\bF_q$ say, then there is  a Frobenius endomorphism of 
$$
S^n_{X^o,x},
$$
which will agree with the corresponding Frobenius endomorphism on 
$$
\pi_1^{\cris}(X^o,x)\ten K.
$$

From now on, assume that $\fX$ is the formal completion of a smooth, proper scheme $X_W/W$, with $\fD$ the formal completion of a divisor $D_W$ of simple normal crossings relative to $W$. Assume moreover that the Frobenius action on 
$$
\H^1_{\cris}(X/W)\ten K
$$
is pure of weight $1$, noting in particular that this will hold if $X$ is projective (by \cite{KM}). 

If we now consider the Gysin sequence (as in Lemma \ref{gdwgt}):
$$
\H^1(X_{\bar{K}},\Q_p) \to \H^1(X^o_{\bar{K}},\Q_p)) \to \H^0(\tilde{D}_{\bar{K}},\Q_p(-1)),
$$
then using \cite{Hop} Theorem 5.6 (the crystalline conjecture), we see that the following sequence is exact:
$$ 
\H^1_{\cris}(X/W)\ten K \to \H^1_{\cris}(X^o/W)\ten K \to \H^0_{\cris}(\tilde{D}/W)\ten K[-1],
$$
so that we have a weight decomposition (due to Frobenius) on $\H^1_{\cris}(X^o/W)\ten K$, splitting the weight filtration.

\begin{theorem} 
The Frobenius action on $\cL(\pi_1^{\cris}(X^o,x), K)$ provides a splitting of weight filtration. Moreover, if  $\H^1_{\cris}(X^o/W)\ten K$ is pure of  weight $a>0$, then there is a Frobenius equivariant map
$$
L((\H^1_{\cris}(X^o/W)\ten K)^{\vee}) \onto \cL(\pi_1^{\cris}(X^o,x), K),
$$
inducing the identity on tangent spaces $\H^1_{\cris}(X^o/W)\ten K$.
\begin{proof}
As for Theorems \ref{finitewgt} and \ref{finitestr}.
\end{proof}
\end{theorem} 

\subsection{The Hodge filtration}\label{hodgefil}

Considering  the commutative diagram on page \pageref{criscompare},  and
given a crystal $\sF$ on $(X/W)_{\cris}$, define the sheaf
$$
\sC^n_{\Hdg}(\sF):= (j^{-1}j_*)^{n+1}\sF.
$$
These give rise to a natural cosimplicial complex of sheaves, which is a resolution of the crystal $\sF$, i.e. is exact for the functor $u_{X^o/W*}$. When $\sF$ is a group-valued sheaf, we can define a Hodge filtration on this cosimplicial complex, by
$$
\Fil^i\sC^n_{\Hdg}(\sF)= \left\{ \begin{matrix} \bigcap_{j_1,\ldots j_{n+1-i}} \ker(\sigma^{j_1}\cdots \sigma^{j_{n+1-i}}:\sC^n_{\Hdg}(\sF) \to \sC^{i-1}_{\Hdg}(\sF)) & n \ge i				\\
						1 & n< i. \end{matrix} \right.
$$
Note that each $\Fil^i\sC^{\bullet}_{\Hdg}(\sF)$ is then a cosimplicial complex of groups, and that this agrees with the usual Hodge filtration on the normal complex, so that if $\sF$ is abelian, this will induce the Hodge filtration on cohomology. 

Now, observe that, if we take $\CC_{\cris}$ as defined on pages  \pageref{crishtpy} and \pageref{logcrishtpy}, then
\begin{eqnarray*}
S^n_{X^o,x}(\g) &=& \CC_{\cris}^n(X^o_{\cris}, \exp(\g\ten \O_{X^o/W}))\\
&=& \Gamma(X^o_{\cris},(w_* w^{-1})^{n+1} \exp(\g\ten \O_{X^o/W}))\\
&=& \Gamma(X^o_{\cris},(v_* v^{-1})^{n+1} (j^{-1}j_*)^{n+1} \exp(\g\ten \O_{X^o/W}))\\
&=& \Gamma(X^o_{\cris},(v_* v^{-1})^{n+1}  \sC^n_{\Hdg}(\exp(\g\ten \O_{X^o/W})))\\
&=& \Gamma(X^o_{\Zar}, u_{X^o/W*} (v_* v^{-1})^{n+1}  \sC^n_{\Hdg}(\exp(\g\ten \O_{X^o/W})))\\
&=& \Gamma(X^o_{\Zar},  (v_* v^{-1})^{n+1} u_{X^o/W*} \sC^n_{\Hdg}(\exp(\g\ten \O_{X^o/W})))\\
&=& \check{\CC}(\cU^o , u_{X^o/W*} \sC^n_{\Hdg}(\exp(\g\ten \O_{X^o/W}))),
\end{eqnarray*}
on which the Hodge filtration is induced from that on $\sC^n_{\Hdg}$.

\section{$p$-adic Hodge Theory}\label{phodge}

In this section we will establish a comparison theorem between the \'etale and crystalline homotopy types. This will imply that the pro-$p$-unipotent fundamental group is a crystalline representation, thus satisfying the hypotheses of Section \ref{pgdl=p}.

Let $V$ be a complete discrete valuation ring, with fraction field $K$, and finite residue field $k$ of characteristic $p$. Denote by $\bar{K}$ (resp. $\bar{k}$) the algebraic closure of $K$ (resp. $k$), and $\bar{V}$ the integral closure of $V$ in $\bar{K}$.  Write $G=\Gal(\bar{K}/K)$. Let $W=W(k)$, the ring of Witt vectors over $k$, and $K_0$ the fraction field of $W$.

$X$ will be a smooth and proper scheme over $V$, and $x$ a $V$-valued point of $X$, with corresponding points $x_{\eta} \in X(K),\, x_s \in X(k)$, and the geometric points \mbox{$\bar{x}_{\eta} \in X(\bar{K}),$} \mbox{$\bar{x}_s \in  X(\bar{k})$.}

We adopt the conventions that
$$
H^*_{\cris}(X):=H^*_{\cris}(X_k/W)\ten_W K_0,
$$
where
$$  
H^*_{\cris}(X_k/W)=\Lim H^*_{\cris}(X_k/W_n),
$$

Following \cite{pi1} Remark \ref{pi1-htpy}, it ought to be true that every cohomology theory with coefficients in a field of characteristic zero should have an associated homotopy type, arising from a standard cosimplicial algebra employed to calculate the cohomology. The SDC arising from this cosimplicial algebra should govern the pro-nilpotent fundamental group.

In the crystalline case, the homotopy type is the complex
$$
\CC^{\bullet}_{\cris}(X, K_0),
$$
defined on page \pageref{crishtpy}.

From \cite{Hop}, there is an isomorphism
$$
\H^*_{\cris}(X)\ten  \Bcr \cong \H^*(X_{\bar{K}},\Q_p)\ten \Bcr,
$$
preserving various structures. It is thus natural to expect that there should be a quasi-isomorphism
$$
\CC^{\bullet}_{\cris}(X)\ten \Bcr \sim  \CC^{\bullet}(X_{\bar{K}},\Q_p)\ten \Bcr
$$
preserving the various structures. This would imply an isomorphism
$$
\pi_1^{\cris}(X_k,x_s)\ten \Bcr \cong \pi_1(X_{\bar{K}},\bar{x}_{\eta})\ten \Bcr.
$$

If $X$ is smooth and projective, then Corollary \ref{crisquad} provided a Frobenius equivariant isomorphism
$$
\cL(\pi_1^{\cris}(X_k,x_s),K_0) \cong  L(\H^1_{\cris}(X)^{\vee})/(\check{\cup}\,\H^2_{\cris}(X)^{\vee})),
$$
where $\check{\cup}$ is dual to the cup product.

This will then provide  a Galois equivariant isomorphism 
$$
\cL(\pi_1(X_{\bar{K}},\bar{x}_{\eta}),\Q_p) \cong L(\H^1(X_{\bar{K}},\Q_p)^{\vee})/(\check{\cup}\H^2(X\ten \bar{K},\Q_p)^{\vee})).
$$

Similarly, there will be a log-crystalline version. All working will be given for  the open case, since the extra generality is no harder to handle.

Let $X^o=X -D$, for $X$ smooth and proper over $V$, and $D$ a divisor with simple normal crossings over $V$. 

\subsection{Comparison of \'etale and intermediate homotopy types}\label{etint}

In this section, we will establish the quasi-isomorphism between the \'etale homotopy type, and the homotopy type corresponding to Faltings' intermediate cohomology (\cite{Hop}). Given a constructible $\Q_p$-sheaf $\sF$ on $X$,  we denote by
$$
\CC^{\bullet}_{\et}(X,\sF) \quad \text{ and } \quad \sC^{\bullet}_{\et}(\sF)
$$
the complexes defined in \cite{pi1} Section \ref{pi1-nrep}.

Given a scheme $X$,  an affine cover $\cU=\{U_{\alpha}\}$  of  $X$, and a Zariski sheaf $\sF$ on $X$, define the associated \v Cech complex
$$
\check{\CC}^n(\cU,\sF):= \Gamma(X, (u_*u^*)^{n+1}\sF),
$$
where $\coprod_{\alpha} U_{\alpha} \xra{u} X$. This has a canonical cosimplicial structure, associated to the adjunction $u^*\dashv u_*$. 
%
%

Given a topological group $\Delta$, and a topological field $F$, define the complex computing group cohomology by
$$
\CC^n_{\gp}(\Delta, F):= \Hom_{\Set,\cts}(\Delta^{n+1},F),
$$
with the usual cosimplicial structure.

In the situation above, take an affine cover $\cU$ of $X$, and denote by $\cU^o$ the corresponding cover of $X^o$. We then have a quasi-isomorphism
\begin{equation}\label{qimetzar}
\CC^{\bullet}_{\et} (X^o\ten \bar{K},\Q_p) \to \check{\CC}^{\bullet}(\cU^o\ten \bar{K},\sC^{\bullet}_{\et}(\Q_p)),
\end{equation}
where, on the right, we mean the diagonal complex associated to the double complex shown, a convention to be adopted throughout this section. 
That this is a quasi-isomorphism follows from the observation that $\sC^n_{\et}(\Q_p)$ is acyclic for Zariski cohomology (since it is for \'etale cohomology), forcing the second spectral sequence (associated to the bigrading)  to degenerate.

As in \cite{Hop} p.46, define the \'etale presheaf of geometric fundamental groups $\Delta$ on $X^0$ by
$$
\Delta(U)=\pi_1(U\ten \bar{K}),
$$
this only being a presheaf, since there is a certain ambiguity about the isomorphisms, due to inner automorphisms. This gives rise to the \'etale presheaf
$$
\sC^{\bullet}_{\gp}(\Delta, \Q_p),
$$
of cosimplicial complexes. 

Since $\sC^{\bullet}_{\gp}(\Delta, \Q_p)$ is a resolution of $\Q_p$, there is a quasi-isomorphism:
\begin{equation}\label{qimetgp}
\check{\CC}^{\bullet}(\cU^o\ten \bar{K},\sC^{\bullet}_{\et}(\Q_p)) \to \check{\CC}^{\bullet}(\cU^o\ten \bar{K},\sC^{\bullet}_{\et}(\sC^{\bullet}_{\gp}(\Delta, \Q_p))),
\end{equation}
where on the right we are taking the diagonal complex associated to a triple cosimplicial complex. 

Now, consider the map
\begin{equation}\label{qimetint}
\check{\CC}^{\bullet}(\cU^o\ten \bar{K},\sC^{\bullet}_{\gp}(\Delta, \Q_p)) \to \check{\CC}^{\bullet}(\cU^o\ten \bar{K},\sC^{\bullet}_{\et}(\sC^{\bullet}_{\gp}(\Delta, \Q_p))).
\end{equation}
In general, this will not be a quasi-isomorphism. However, it will be whenever the $U_{\alpha}$ are such that the $U_{\alpha}^o\ten \bar{K}$ are $K(\pi,1)$ spaces. 

Recall (for instance, from \cite{FaH} Ch.II \S 2) the definition of a $K(\pi,1)$ space.  Suppose $Z$ is a smooth irreducible scheme over an algebraically closed field of characteristic zero, with fundamental group $\pi$. If $F$ is a finite abelian group with a continuous $\pi$-action, corresponding to a locally constant sheaf $\mathbf{F}$. There is a natural transformation
$$
\H^*_{\gp}(\pi,F) \to \H^*_{\et}(X,\mathbf{F})
$$
(corresponding, in fact to the morphism (\ref{qimetint}) of complexes above). $Z$ is said to be a $K(\pi,1)$ if this map is an isomorphism for all such $F$.

By \cite{FaH} Ch.II Lemma 2.1,  there exist affine covers $\cU$ of $X$ for which each $U_{\alpha}^o\ten \bar{K}$ is a $K(\pi,1)$. It then follows that the map (\ref{qimetint}) is a quasi-isomorphism for such covers. 

We thus obtain a complex computing intermediate cohomology by
$$
 \check{\CC}^{\bullet}(\cU^o\ten \bar{K},\sC^{\bullet}_{\gp}(\Delta, \Q_p)),
$$
for $\cU^o\ten\bar{K}$ a cover by open affine $K(\pi,1)$'s, which we have hence shown to be quasi-isomorphic (as a cosimplicial algebra) to  
$$
\CC^{\bullet}_{\et} (X^o\ten \bar{K},\Q_p).
$$
 Note that since all maps were canonical, this quasi-isomorphism is $\Gal(\bar{K}/K)$-invariant.

\subsection{Comparison of crystalline and intermediate homotopy types}\label{crisint}

Recall from Section \ref{etint} that the intermediate homotopy type is
$$
\CC^{\bullet}_{\mathrm{int}}(X^o,\Q_p) =\check{\CC}^{\bullet}(\cU^o\ten \bar{K},\sC^{\bullet}_{\gp}(\Delta, \Q_p)),
$$
while, from Section \ref{hodgefil}, the crystalline homotopy type is
$$
\CC^{\bullet}_{\cris}(X_k^o,K_0) \cong \check{\CC}(\cU^o_k , u_{X^o_k/W*} \sC^n_{\Hdg}(\O_{X^o_k/W}\ten K_0)).
$$

As in \cite{Hop} p.62,  we may use the crystal $\sC^{\bullet}_{\Hdg}(\O_{X^o/W_n})$ to  define sheaves
$$
\sC^{\bullet}_{\Hdg}(\O_{X_k^o/W_n})\ten_B^+/(I^{[n]}+p^nB^+)
$$
on the site $\cX^o$ defined in \cite{Hop}, by associating to each $U=\Spec R$, the module 
$$
u_{X_k^o/W_n*}\sC^{\bullet}_{\Hdg}(\O_{X_k^o/W_n})(R)\ten_R B^+(\hat{R})/(I^{[n]}+p^nB^+).
$$

This then gives us a morphism
\begin{eqnarray*}
&&
\check{\CC}(\cU^o_k, u_{X^o/W_n*} \sC^n_{\Hdg}(\O_{X^o_k/W_n}))\ten B^+(V)/(I^{[n]}+p^nB^+)\\ 
&\to& \check{\CC}^{\bullet}(\cU^o\ten \bar{K},\sC^{\bullet}_{\gp}(\Delta,u_{X^o_k/W_n*}\sC^{\bullet}_{\Hdg}(\O_{X^o_k/W})\ten_B^+/(I^{[n]}+p^nB^+) )),
\end{eqnarray*}
which, by taking inverse limits and tensoring (over $B^+(V)$) with $B(V)$, gives a morphism
$$
\CC^{\bullet}_{\cris}(X^o,K_0)\ten_{K_0}B(V) \to  \Lim \CC^{\bullet}_{\mathrm{int}}(X^o,u_{X^o_k/W_n*}\sC^{\bullet}_{\Hdg}(\O_{X^o_k/W})\ten_B^+/(I^{[n]}+p^nB^+))\ten B(V)
$$ 
which by \cite{Hop}  is a filtered quasi-isomorphism, as is the  natural morphism
$$
\CC^{\bullet}_{\mathrm{int}}(X^o,\Q_p)\ten_{\Q_p}B(V) \to  \Lim \CC^{\bullet}_{\mathrm{int}}(X^o,u_{X^o_k/W_n*}\sC^{\bullet}_{\Hdg}(\O_{X^o_k/W})\ten_B^+/(I^{[n]}+p^nB^+))\ten B(V).
$$

We therefore have a quasi-isomorphism
$$
\CC^{\bullet}_{\cris}(X_k^o,K_0)\ten_{K_0} B(V) \sim \CC^{\bullet}_{\mathrm{int}}(X_{\bar{K}^o},\Q_p)\ten_{\Q_p} B(V),
$$
respecting the Hodge filtration, and the Frobenius and Galois actions.

\subsection{Consequences for fundamental groups}

Combining the two previous sections, we have a quasi-isomorphism
$$
\CC^{\bullet}_{\cris}(X_k^o,K_0)\ten_{K_0} B(V) \sim \CC^{\bullet}_{\et}(X_{\bar{K}}^o ,\Q_p)\ten_{\Q_p} B(V),
$$
preserving the various structures.

If $x \in X^o(V)$, let $x_{\eta} \in X^o(K)$ and $x_s \in X^o(k)$ be the corresponding closed points, with $\bar{x}_{\eta} \in X^o(\bar{K})$ the geometric point.

Recall that
$$
S^{\bullet}_{\et}(\g):= \exp(\CC^{\bullet}_{\et}(X_{\bar{K}}^o,\bar{x}_{\eta};\Q_p)\ten_{\Q_p}\g)
$$
is the SDC over $\cN_{\Q_p}$ giving rise to the Lie algebra
$$
\cL^{\et}:=\cL(\pi_1^{\et}(X_{\bar{K}}^o,\bar{x}_{\eta}),\Q_p),
$$
while the SDC 
$$
S^{\bullet}_{\cris}(\g):= \exp(\CC^{\bullet}_{\cris}(X_k^o,x_s;K_0)\ten_{K_0}\g)
$$
over $\cN_{K_0}$ gives rise to the Lie algebra
$$
\cL^{\cris}:=\cL(\pi_1^{\cris}(X_k^o,x_s),K_0).
$$

We wish to compare these two.

\begin{definition}\label{ninfty} Given a field $\kappa$ of characteristic zero, define the category $\cN_{\kappa}^{\infty}$ to consist of nilpotent Lie algebras over $\kappa$, not necessarily finite dimensional. Note that any functor $F$ on $\cN_{\kappa}$ which commutes with filtered direct limits can be extended uniquely to $\cN_{\kappa}^{\infty}$ by setting
$$
F(\g)= \varinjlim F(\g_{i}),
$$
where the $\g_i$ run over all finitely generated (and hence finite dimensional) sub-Lie algebras of $\g$.
\end{definition}

 Recall that, on $\cN_{\Q_p}$, we have functorial isomorphisms
$$
\Hom_{\widehat{\cN}_{\Q_p}}(\cL^{\et},\g) \cong \ddef_{S_{\et}}(\g).
$$
Now, the functor $S^n_{\et}$ extends naturally to $\cN_{\Q_p}^{\infty}$, and commutes with direct limits since tensor products do. The functor $\Hom(\cL^{\et},-)$ also extends naturally to $\cN_{\Q_p}^{\infty}$, and commutes with direct limits since the Lie algebra is pro-finitely presented.

Note that there is a canonical map
$$
\cN_{K_0}^{\infty} \to \cN_{\Q_p}^{\infty},
$$
given by the forgetful functor, so that the isomorphism of functors above can be used to give an isomorphism of functors on $\cN_{K_0}^{\infty}$, namely
$$
\Hom_{\widehat{\cN}_{K_0}^{\infty}}(\cL^{\et}\ten_{\Q_p}K_0,\g) \cong \ddef_{S_{\et}\ten_{\Q_p} K_0}(\g),
$$
where
$$
S^n_{\et}\ten_{\Q_p} K_0(\g):=\exp(\CC^{\bullet}_{\et}(X^o_{\bar{K}},\bar{x}_{\eta};\Q_p)\ten_{\Q_p}\g).
$$

Similarly, we have an isomorphism
$$
\Hom_{\widehat{\cN}_{K_0}^{\infty}}(\cL^{\cris},\g) \cong \ddef_{S_{\cris}}(\g).
$$

\begin{theorem}
There is a canonical isomorphism
$$
\Hom_{\widehat{\cN}_{K_0}^{\infty}}(\cL^{\cris},\g\ten_{K_0} B) \cong \Hom_{\widehat{\cN}_{K_0}^{\infty}}(\cL^{\et}\ten_{\Q_p}K_0,\g\ten_{K_0} B),
$$
functorial in $\g \in \cN_{K_0}^{\infty}$, and respecting the Galois and Frobenius actions. In fact, this is functorial not only in $\g$, but also in $\g\ten_{K_0} B$, meaning that the isomorphism respects $B$-linear Lie algebra maps $\g \ten_{K_0} B \to \fh \ten_{K_0} B$.
\begin{proof}
The isomorphisms established above imply that it suffices to show 
$$
\ddef_{S_{\cris}}(\g\ten_{K_0} B) \cong \ddef_{S_{\et}\ten_{\Q_p} K_0}(\g \ten_{K_0} B),
$$
functorially in $\g\ten_{K_0} B$.

Now, observe that
$$
\ddef_{S_{\cris}}(\g\ten_{K_0} B) \cong \ddef_{S_{\cris}\ten_{K_0}B}(\g),
$$
where
$$
S^n_{\cris}\ten_{K_0}B(\g):=\exp((\CC^{\bullet}_{\cris}(X_k^o,x_s;K_0)\ten_{K_0}B)\ten_{K_0}\g).
$$
Similarly
$$
\ddef_{S_{\et}\ten_{\Q_p} K_0}(\g\ten_{K_0} B) \cong \ddef_{S_{\et}\ten_{\Q_p}B}(\g),
$$
where
$$
S^n_{\et}\ten_{\Q_p} B(\g):=\exp((\CC^{\bullet}_{\et}(X^o_{\bar{K}},\bar{x}_{\eta};\Q_p)\ten_{\Q_p}B)\ten_{K_0}  \g).
$$

Now the quasi-isomorphism 
$$
\CC^{\bullet}_{\cris}(X_k^o,x_s;K_0)\ten_{K_0} B \sim \CC^{\bullet}_{\et}(X_{\bar{K}}^o,\bar{x}_{\eta}; \Q_p)\ten_{\Q_p}
$$
established in the previous sections provides us with an isomorphism of the associated deformation functors
$$
\ddef_{S_{\cris}\ten_{K_0}B}(\g) \cong \ddef_{S_{\et}\ten_{\Q_p}B}(\g),
$$
for  $\g \in \cN_{K_0}$. Now, since tensor products commute with direct limits, both functors do, giving us isomorphisms
$$
\ddef_{S_{\cris}\ten_{K_0}B}(\g) \cong \ddef_{S_{\et}\ten_{\Q_p}B}(\g),
$$
 for all $\g \in \cN_{K_0}^{\infty}$, as required. 
\end{proof}
\end{theorem}

\begin{corollary}\label{compthm} There is a Galois and Frobenius equivariant isomorphism
$$
\cL(\pi_1^{\cris}(X_k^o,x_s),K_0)\ten_{K_0}B \cong \cL(\pi_1^{\et}(X_{\bar{K}}^o,\bar{x}_{\eta}),\Q_p)\ten_{\Q_p} B.
$$
\begin{proof} From the theorem, we have an isomorphism
$$
\Hom_{\widehat{\cN}_{K_0}^{\infty}}(\cL^{\cris},\g\ten_{K_0} B) \xra{\theta} \Hom_{\widehat{\cN}_{K_0}^{\infty}}(\cL^{\et}\ten_{\Q_p}K_0,\g\ten_{K_0} B),
$$
or equivalently an isomorphism
$$
\Hom_{B}(\cL^{\cris}\ten_{K_0}B,\g\ten_{K_0} B) \xra{\theta} \Hom_{B}(\cL^{\et}\ten_{\Q_p}B,\g\ten_{K_0} B),
$$
functorial not only in $\g$, but also in $\g\ten_{K_0}B$.
Define
$$
\alpha=\theta(\id): \cL^{\et}\ten_{\Q_p}B \to \cL^{\cris} \ten_{K_0} B,
$$
and define
$$
\beta=\theta^{-1} (\id):\cL^{\cris}\ten_{K_0}B \to \cL^{\et}\ten_{\Q_p}B.
$$

Now, for $f: \g\ten_{K_0} B \to \fh\ten_{K_0} B$, and $a:\cL^{\cris}\ten_{K_0}B \to \g\ten_{K_0} B$, we must have
$$
\theta(f\circ a)= f\circ \theta(a),
$$
by functoriality, and similarly for $\theta^{-1}$.

Therefore
\begin{eqnarray*}
\id =& \theta(\beta)=\theta(\beta \circ \id)&=\beta \circ \alpha\\
\id =& \theta^{-1}(\alpha)=\theta^{-1}(\alpha \circ \id)&=\alpha \circ \beta.
\end{eqnarray*}
\end{proof}
\end{corollary}

We now have the following corollary, showing that in this case the hypotheses of Section \ref{pgdl=p} are satisfied:

\begin{corollary}\label{crys}
If $X^o=X -D$, for $X$ smooth and proper over $V$, and $D$ a divisor with simple normal crossings over $V$, then  $\cL(\pi_1^{\et}(X^o_{\bar{K}},\bar{x}_{\eta}),\Q_p)$ is a crystalline Galois representation.
\end{corollary}

\begin{remark} 
 A generalisation of this comparison result to relative Malcev completions has since been independently obtained by Olsson (\cite{olssonhodge} Theorem 1.9). 
\end{remark}

\newpage
\bibliographystyle{alpha}
\addcontentsline{toc}{section}{Bibliography}
\bibliography{references.bib}
\end{document}